\documentclass[12pt]{article}

\usepackage{amssymb,amsmath,amsthm}
\usepackage{epsfig}

\makeatletter

\@addtoreset{equation}{section}
\makeatother

\catcode`\@=11
\def\qed{{\unskip\nobreak\hfil\penalty50
          \hskip2em\hbox{}\nobreak\hfil\mbox{\rule{1ex}{1ex}}
          \parfillskip=0pt
   \finalhyphendemerits=0\par\medskip}}
\catcode`\@=12

\def\ocirc#1{\ifmmode\setbox0=\hbox{$#1$}\dimen0=\ht0 \advance\dimen0
  by1pt\rlap{\hbox to\wd0{\hss\raise\dimen0
  \hbox{\hskip.2em$\scriptscriptstyle\circ$}\hss}}#1\else {\accent"17 #1}\fi}
\def\cprime{$'$}

\def\R{\mathbf{R}}
\def\ds{\displaystyle}

\def\div{{\rm div}}

\def\refe#1{(\ref{#1})}

\def\Lip{\mathrm{Lip}}

\def\bbarint#1#2{{\int \!\!\!\!\!\!-}_{\!\!\!\! #1}^{#2}}

\newtheorem{theorem}{Theorem}
\newtheorem{lemma}{Lemma}[section]
\newtheorem{corollary}{Corollary}[section]

\newtheorem{proposition}{Proposition}[section]

\theoremstyle{remark}
\newtheorem{remark}{Remark}[section]

\DeclareMathOperator{\supp}{supp}

\begin{document}

\title{Error estimate for the Finite Volume scheme applied to the linear advection
equation}
\author{Beno{\^\i}t Merlet and Julien Vovelle}
\date{}
\maketitle


{\bf Summary. }  We study the convergence of a Finite Volume scheme
  for the linear advection equation with a  Lipschitz divergence-free
  speed in $\R^d$. We prove a $h^{1/2}$-error estimate in the $L^\infty(0,t;L^1)$-norm for $BV$ data. This result was expected from numerical experiments and is optimal. 

\medskip


\medskip


\section{Introduction}\label{s1}

The Finite Volume method is well adapted to the computation of the solution of pdes which are 
conservation (or balance) laws, for the reason that it respects the
property of conservation of the pde under study. 
The mathematical analysis of the application 
of the Finite Volume method to hyperbolic first-order conservation laws can be dated from the mid
sixties (see \cite{TihonovSamarskii62} for example). Concerning the specific problem of the estimate 
of the rate of convergence of the method, the first result is due to 
Kuztnetsov~\cite{Kuznetsov76}, who proves that this rate of convergence in $L^\infty(0,t;L^1)$ 
is of order $h^{1/2}$, where $h$ is the size of the mesh, provided that the initial data is in $BV$ and
that the mesh is a structured cartesian grid. Ever since, several studies and results have come to 
supplement the error estimate of Kuznetsov. Before describing them, let us emphasize two points:
\medskip

{\it 1.} The analysis of the speed of convergence of the Finite Volume method is distinct from the
analysis of the order of the method. In the analysis of the speed of convergence of the method, general data (e.g. $BV$ data) are 
considered. Indeed, here, the problem is to show that the Finite Volume method behaves well regarding
the approximation of the continuous evolution problem in all his features (in particular the creation and the transport of discontinuities).
On the other hand, in the analysis of the order of the method, 
restrictions on the regularity of the data are of no importance: see the recent work of Bouche, Ghidaglia, 
Pascal~\cite{BoucheGhidagliaPascal05} on that purpose.
\medskip

{\it 2.} If, numerically, the speed of convergence of the Finite Volume method applied to first-order
conservation laws is observed to be (at least) of order $h^{1/2}$ in the $L^\infty(0,t;L^1)$-norm,
whether the mesh is structured or is not, the preexisting theoretical and rigorous proofs of this result appears to be 
strongly related to the structure of the mesh.
\medskip

Indeed, in the case where the mesh is unstructured (what we call
a structured mesh is a cartesian mesh with identical cells but this can be slightly relaxed~\cite{CockburnGremaudYang98}), 
the result of Kuznetsov has been extended, but to
the price of a fall in the order of the error estimate. Indeed, $h^{1/4}$ error estimate in the 
$L^\infty(0,t;L^1)$-norm for
the Finite Volume method applied to hyperbolic conservation laws on unstructured meshes has been proved
by Cockburn, Coquel, Lefloch \cite{CockburnCoquelLefloch94}, Vila \cite{Vila94} and Eymard, Gallou\"et,
Herbin \cite{EGH00} for the Cauchy Problem ($h^{1/6}$ error estimate for the Cauchy-Dirichlet Problem
\cite{OhlbergerVovelle04}). We emphasize the fact that numerical tests give an order $h^{1/2}$ for structured as
well as unstructured meshes; still, concerning these latter, numerical analysis did not manage to give the
rigorous proof of the order $h^{1/2}$: there is an upper limit at the order $h^{1/4}$. 

\medskip

In this paper, we consider the case where the conservation law is
{\it linear}. More precisely, we consider the linear advection problem
with a Lipschitz divergence-free speed. 

In the case of $H^1$-data, optimal (with respect to the order of the
error estimate) results already exist. Under a strict CFL condition,
Despr\'es~\cite{Despres04nl} shows an $h^{1/2}$-error estimate in the $L^\infty_t
L^2_x$-norm for the upwind Finite Volume method applied to the linear
advection equation with constant speed. (The author deals with the particular case of 2D triangular
meshes but extensions to higher dimension and general polyhedral
meshes are harmless). His technique is based on the study of the
dissipation of the consistency error by the scheme and the adaptation
of the Lax Theorem (see also~\cite{Despres04lax}). Roughly speaking, the consistency error
created at time $n\delta t$ is of order $1$ but of order
$1/\sqrt{q+1}$ after $q$ time steps of the scheme.  

On the other hand, we refer to the work of Vila and Villedieu
\cite{VilaVilledieu03}, who prove, again under a strict CFL condition, an $h^{1/2}$-error estimate in the $L^2_{loc}$ space-time norm for the
approximation of Friedrichs hyperbolic systems with $H^1$ data by energy estimates (they consider explicit in time Finite Volume schemes, for
implicit schemes in the case of scalar advection equation, see, as they underline it,
the result of Johnson and Pitk{\"a}ranta \cite{JohnsonPitkaranta86} who show an $h^{1/2}$-error estimate in the $L^2$ space-time norm for $H^1$ data). 
\smallskip

For an initial data $u_0\in BV(\R^d)$, we prove in Theorem~\ref{thmL1}
the expected $h^{1/2}$-error estimate in the
$L^\infty_t L^1_x$ norm, under a strict CFL
condition in the general case and, under a sharp CFL condition if the
speed is independent of $t$.  This result is optimal \cite{TangTeng95,Sabac97}
and in the context of unstructured meshes and $BV$-data, this is the first optimal result.

In our proof of Theorem~\ref{thmL1}, we first show that it is
sufficient to consider initial data which are characteristic
functions. Then we use the same ingredients as Vila and Villedieu
in~\cite{VilaVilledieu03}. The approximate solution satisfies
the weak formulation of the problem up to an error term corresponding
to the consistency error of the scheme. Roughly speaking, this error is bounded via
the Cauchy-Schwarz inequality by a $H^1(\R^d\times[0,t])$-like semi-norm
of the discrete solution denoted $E_h(u_0,t)$. A large error means
that the scheme is very dissipative and consequently, that the approximate
solution is smooth, so $E_h(u_0,t)$ is small. The $h^{1/2}$-estimate
follows from the equilibrium between these contradictory constraints.

\medskip

The paper is divided into six parts: first we continue this
introduction by describing the linear advection problem, the Finite
Volume scheme, our results and we give the main lines of the proof. In
Section~\ref{s2} we recall some classical results on the Finite Volume
scheme. In Section~\ref{s3} we introduce the weak formulation
satisfied by the approximate solution. Then we
reduce the study to initial data which are characteristic functions and we build the test function used in the weak
formulation. In Section~\ref{s4}, we prove
some Energy estimates. The first and second part of
Theorem~\ref{thmL1} are proved in Sections~\ref{s41},~\ref{s42} respectively.

\subsection*{Notations}
If $(X,\mu)$ is a measurable set with finite (positive) measure and $\phi\in L^1(X)$, we denote the mean of $\phi$ over $X$ by
$$
\ds\bbarint{X}{} \phi d\mu:=\frac{1}{\mu(X)}\int_X \phi d\mu.
$$
\smallskip

If $X$ is a set, $\mathbf{1}_X(x)=1$ if $x\in X$, $0$ if $x\notin X$. 
\smallskip

The set of functions with bounded variation in $\R^d$ is the set of $L^1$ functions with bounded Radon measures as derivatives:
$$
\ds BV(\R^d):=\left\{u\in L^1(\R^d);\sup_{\varphi\in{\cal C}^\infty_c(\R^d,\R^d),\|\varphi\|_\infty\leq 1}\left|\int_U u(x){\rm div}\varphi(x) dx\right|
<+\infty\right\}
$$
where $||\varphi||_\infty:=||\sqrt{\varphi_1^2+\cdots+\varphi_d^2}||_{L^\infty(\R^d)}$ for $\varphi\in{\cal C}^\infty_c(\R^d,\R^d)$ and ${\rm div}$
is the divergence operator. The total variation of $u\in BV(\R^d)$ is
given by 
\begin{eqnarray*}
\|u\|_{TV} &:=& \sup_{\varphi\in{\cal C}^\infty_c(\R^d,\R^d),\|\varphi\|_\infty\leq 1}\left|\int_{\R^d} u(x){\rm div}\varphi(x) dx\right|.
\end{eqnarray*}  

\medskip



\subsection{The linear advection equation}\label{secPB}

We consider the linear advection problem in $\R^d$:

\begin{equation}\label{PB}
\left\{\begin{array}{l l}
u_t+{\rm div}(V u)=0, & x\in\R^d, t\in\R_+,\\
\\
u(x,0)=u_0(x), & x\in\Omega, \\
\end{array}\right.
\end{equation}
where we suppose that $V \in W^{1,\infty}(\R^d \times \R_+,\R^d)$ satisfies  
\begin{equation}
\label{divergence}
\div \, V(\cdot,t)=0 \qquad \forall t\in \R_+.
\end{equation}
The problem \refe{PB} has a solution for $u_0\in L^1_{loc}(\R^d)$; for the purpose of the error estimates, 
we will consider initial data in $BV(\R^d)$.

\begin{theorem}\label{thsolPB} For every $u_0\in L^1_{loc}(\R^d)$, the problem~\refe{PB} admits a unique weak solution $u$, in the sense
that $u\in L^1_{loc}(\R_+\times \R^d)$ and: for every $\phi\in{\cal
  C}_c^\infty({\R^d}\times \R_+)$,
\begin{equation}
\label{formfaible}
\ds\int_{\R_+}\int_{\R^d} u(\phi_t+V\cdot\nabla\phi) dx dt +\int_{\R^d} u_0\phi(x,0) dx =0.
\end{equation}
Moreover, we have 
\begin{equation}
\label{formulexacte}
u(x,t)=u_0(X(x,t)), \qquad \forall (x,t) \in \R^d \times \R_+,
\end{equation}
where $X\in C^1(\R^d\times \R_+, \R^d)$ is such that, for every $t\geq0$, $X(\cdot,t):\R^d\rightarrow \R^d$ is one to one and onto.

 Let $Y(\cdot,t):=X(\cdot,t)^{-1}$, then for every $T\geq 0$,  $X-Id_{\R^d_x}$ and  $Y-Id_{\R^d_x}$ belong to  $W^{1,\infty}(\R^d \times [0,T])$. 
More precisely, there exists $C_0 \geq 1$ only depending on $T$, $\|V\|_{W^{1,\infty}}$ and $d$, such that, $\forall x\in \R^d,\,\forall t\in[0,T]$,
\begin{eqnarray} 
\label{estimationX}
\|X(x,t)-x\|\leq C_0\, t ,& \qquad \|\nabla X(x,t)\| \leq C_0,&\qquad \|\partial_tX(x,t)\| \leq C_0,\\
\|Y(x,t)-x\|\leq C_0\,t ,& \qquad \|\nabla Y(x,t)\| \leq C_0,&\qquad \|\partial_t Y(x,t)\| \leq C_0.
\label{estimationY}
\end{eqnarray}

Finally, for every $t\geq 0$, $X(\cdot,t)$ preserves the Lebesgue measure $\lambda$ on $\R^d$, i.e: 
\begin{eqnarray}
\label{conservemesure}
\lambda(X(E,t)) &=&\lambda(E), \qquad \mbox{ for every Borel subset }
E \mbox{ of } \R^d.
\end{eqnarray}
\end{theorem}

{\bf Proof of Theorem~\ref{thsolPB}} All the results cited in the Theorem follow from the characteristic formula:
$u(x,t)=u_0(X(x,t))$ where $X(\cdot,t)=Y(\cdot,t)^{-1}$ and $Y$ solves
the Cauchy Problem $\partial_t Y(x,t)= V(Y(x,t),t)$,
$Y(x,0)=x$. To prove the estimates on $X$, we notice that
$X(x,t)=Z(0;x,t)$ where $Z(\tau;x,t)$ denotes the solution of the Cauchy Problem
\begin{equation*}\label{odechi}
\left\{\begin{array}{l l}
\ds\frac{dZ}{d\tau}(\tau;x,t)=V(Z(\tau;x,t),\tau), & \tau\in[0,t],\\
\\
Z(t;x,t)=x.
\end{array}\right.
\end{equation*}
We only give the sketch of the proof: global existences for $Y$ and
$Z$ and the desired estimates follow from the Cauchy-Lipschitz
Theorem and the Gronwall Lemma. Besides, by \refe{divergence}, the flow preserves the Lebesgue measure on $\R^d$. Existence for \refe{PB} follows from the characteristic formula and, by an argument of duality, uniqueness also. \qed

If $u_0\in L^1_{loc}(\R^d)$, the solution $u$ belongs to
$C(\R_+,L^1_{loc}(\R^d))$. If the derivatives of $u_0$ are bounded
Radon measures, we have a more precise result. Namely, using the estimates above, we obtain
\begin{corollary}
\label{regularite}
Let $T\geq 0$, there exists a constant $C_0\geq 1$ depending on
$\|V\|_{W^{1,\infty}}$, $T$ and $d$ such that, if $u_0\in BV(\R^d)$:
\begin{eqnarray*}
\|u(\cdot,s)-u(\cdot,t)\|_{L^1}&\leq&C_0
 \|u_0\|_{TV} |s-t|
  ,\qquad  0\leq s,t \leq T,\\
\|u(\cdot,t)\|_{TV}&\leq& C_0
 \|u_0\|_{TV},\qquad\qquad  0\leq t\leq T.
\end{eqnarray*}   
\end{corollary}
The following result is a direct consequence of the conservative
property~\eqref{conservemesure}. The case $f(v)=v^2$ will be crucial in the
proof of the main result of this paper.
\begin{corollary}
\label{conservenormes}
Let $f :\R\rightarrow\R$ be a measurable function and $u_0 \in
L^1_{loc}(\R^d)$ such that $f\circ u_0 \in L^1(\R^d)$. Then the quantity
$\ds \int_{\R^d} f(u(\cdot,t))$ is constant.  
\end{corollary}



\subsection{Finite Volume scheme}\label{secFVscheme}

The Finite Volume scheme which approximates \refe{PB} is defined on a mesh
${\cal T}$ which is a family of closed connected polygonal subsets with disjoint
interiors covering $\R^d$; the time half-line is meshed by regular cells of size $\delta t>0$. We also suppose that the partition ${\cal T}$ satisfies the following properties:
the common interface of two control volumes is included in an hyperplane 
of $\R^d$ and
\begin{equation}\label{eqgentilmesh}
\mbox{there exists } \alpha>0 \mbox{ such that }\left\{
\begin{array}{l}
\alpha h^d \leq |K|, \\
|\partial K|\leq {\alpha}^{-1}\, h^{d-1}\,,\; \forall K\in{\cal T},
\end{array}
\right.
\end{equation}
where $h$ is the size of the mesh: $h:=\sup\{diam(K),\, K\in{\cal T}\}$, $|K|$ is the $d$-dimensional Lebesgue measure of $K$ and $|\partial K|$ is the $(d-1)$-dimensional Lebesgue measure of $\partial K$. 
If $K$ and $L$ are two control volumes having a common edge,
we say that $L$ is a neighbor of $K$ and denote (quite abusively) $L\in\partial K$. We also
denote $K\!\!\shortmid\!L$ the common edge and $\mathbf{n}_{KL}$ the
unit normal to $K\!\!\shortmid\!L$ pointing outward $K$. We set
$$
V_{KL}^n:=\int_{K\shortmid L}\bbarint{n\delta t}{(n+1)\delta t}V\cdot \mathbf{n}_{KL}
$$
We denote by $K_n:=K\times[n\delta t,(n+1)\delta t)$
a generic space-time cell, by ${\cal M}:={\cal T}\times\mathbb{N}$ the
space-time mesh and by $\partial K^-_n$ the set 
$$
\partial K_n^-:=\{L\in\partial K,V^n_{KL}<0\}.
$$
In the same way we set $K\!\!\shortmid\!L_n := K\!\!\shortmid\!L \times [n\delta t,(n+1)\delta t)$. 
\medskip

We will assume that the so called Courant-Friedrich-Levy condition is satisfied:
\begin{eqnarray}
\label{CFL}
\sum_{L \in \partial K_n^-} \delta t |V_{KL}^n| & \leq &(1-\xi)
|K|,\qquad \forall K_n \in \mathcal{M},\qquad\mbox{ for some }\, \xi \in [0,1).
\end{eqnarray}

\begin{remark} Under condition \refe{eqgentilmesh}, the CFL condition \refe{CFL} holds as soon as
$$
\|V\|_\infty \delta t\leq \ds(1-\xi)\alpha^{-2}\,h.
$$
\end{remark}
Notice that if $\|V\|_\infty$ is small, we may choose a large time
step $\delta t$. In order to avoid the occurrence of terms with factor ${\delta t}\,{h}^{-1}$ in our estimates, we add to~\eqref{CFL}
the following condition: there exists $c_0\geq 0$ such that 
\begin{eqnarray}
\label{CFL2}
\delta t \leq c_0 h.
\end{eqnarray}

The Finite Volume scheme with explicit time-discretization is defined by the following set of equations:
\begin{eqnarray}
u^0_K=\bbarint{K}{}\, u_0(x)\, dx\,,&\; \forall K\in {\cal T},\label{eqscheme0}\\
\nonumber\\
\ds\frac{u^{n+1}_K-u^n_K}{\delta t}+\ds\frac{1}{|K|}\sum_{L\in\partial K^-_n} V_{KL}^n(u^n_L-u^n_K)=0,&\qquad \forall K_n\in{\cal M}.\label{eqscheme}
\end{eqnarray}

We then denote by $u_h$ the approximate solution of \refe{PB} defined by the Finite Volume scheme:
\begin{equation}\label{equh}
u_h(x,t)=u^n_K,\quad\forall (x,t)\in K_n.
\end{equation}

\subsection*{Main results}
From now on, we assume that $V \in W^{1,\infty}(\R^d \times \R_+,\R^d)$
satisfies~\refe{divergence}. We fix a mesh $\mathcal{T}$ of mesh-size
$h > 0$ satisfying the uniformity conditions~\eqref{eqgentilmesh}. We also fix a
time step $\delta t$ such that the CFL conditions
\eqref{CFL}-\eqref{CFL2} hold.

 We fix two times $0\leq t \leq T$ and we assume that $t=(N+1)\delta
 t$ for some $N \in \mathbb{N}$ (so that
$u_h(\cdot,t) =\sum_{K \in {\cal T}} u_K^{N+1} \mathbf{1}_{K}$). In order to avoid the occurrences of ${h}\,{t}^{-1}$ terms in the estimates we assume that  
\begin{eqnarray}
\label{c1}
h&\leq& c_1 t,\qquad \mbox{ for some }  c_1\geq 0.
\end{eqnarray}
In the sequel, given an initial condition $u_0 \in BV(\R^d)$, the function $u\in C(\R_+,L^1(\R^d))$
denotes the exact solution to~\refe{PB} and $u_h$ its numerical approximation obtained by the upwind Finite Volume method
\refe{eqscheme0}-\refe{eqscheme}-\refe{equh}.

In the results and in the proofs, $C_0\geq 1$ is the constant introduced in
Theorem~\ref{thsolPB} and Corollary~\ref{regularite}; this constant
only depends on $T$, $\|V\|_{W^{1,\infty}}$ and $d$. The letter $C$ denotes various constants which are non
decreasing functions of $\alpha^{-1}$, $c_0$, $c_1$,
$\|V\|_{W^{1,\infty}}$ and $d$ but do not depend on
$h$, $\delta t$, $\xi$, $u_0$, $t$ or  $T$. 

Finally, before stating the main results of the paper, we set 
$${\cal M}_N:=\{K_n\in {\cal M}\; :\; 0 \leq n\leq N\},$$ 
and we introduce the quantities 
\begin{eqnarray*}
E_h(u_0,t)&:=&\sum_{K_n \in \mathcal{M}_N} |K| |u_K^{n+1}-u_K^n|^2 +
\sum_{K_n \in \mathcal{M}_N} \sum_{L \in \partial K_n^-} \delta t
|V_{KL}^n|\, |u_L^n-u_K^n|^2,\\
\mathcal{E}_h(u_0,t)&:=& \|u_h(\cdot,0)\|_{L^2}^2  -\|u_h(\cdot,t)\|_{L^2}^2.
\end{eqnarray*}
\begin{remark} 
Notice that $h^{-2} E_h(u_0,t)$ is a discrete version of the quantity 
$$
\ds \int_{\R^d\times [0,t]} \left(|\partial_t u_h(x,s)|^2 +
\sum_{i=1}^d |V_i(x,s)| |\partial_{x_i} u_h(x,s)|^2 \right)dxds.
$$ 
\end{remark}
We have the following error estimate:

\begin{theorem}
Let $u_0$ in $BV(\R^d)$. Under a strict CFL condition (\refe{CFL} with $\xi >0$), we have :
\begin{eqnarray}
\|u(\cdot,t)-u_h(\cdot,t)\|_{L^1}& \leq &CC_0^{d+2}\xi_*^{-1}\|u_0\|_{TV} (t^{1/2}h^{1/2} + \xi_*^{1/2}\,th), \label{thmL1_1}
\end{eqnarray}
with $\xi_*=\xi$.

Moreover, if $V$ does not depend on the time variable, the estimate is valid with
$\xi_*=1$ uniformly in $\xi\in [0,1)$. 
\label{thmL1}
\end{theorem}

This result deserves some comments.
\begin{remark}
If the speed does not depend on $x$: $V(x,s):=\mathbf{V}(s)$, we have
$X(x,s)=x-\int_0^s\mathbf{V}(r)dr$ and we can choose
$C_0=d^{1/2}+\|\mathbf{V}\|_\infty$ (which does not depend on $T$) in Theorem~\ref{thsolPB} and Corollary~\ref{regularite}. 
\end{remark}

\begin{remark}
In the proof of the first part of Theorem~\ref{thmL1}, the dependency on $\xi$ only appears in Lemma~\ref{lemdeuxcasL2}. 
\end{remark}

\begin{remark}[On the $BV$-norm of the approximate solution]
\label{rkweakBV} 
The possible irregularity of the mesh is an obstacle to the existence of uniform (with respect to $h$) bound on the $L^1_t BV_x$-norm 
$$
\sum_{K_n \in \mathcal{M}_N}\sum_{L\in\partial K_n^-}\delta
t |K\!\!\shortmid\!L| |u^n_L-u^n_K|
$$
of the approximate solution $u_h$. A counter-example to such a result has been given by B.~Despr\'es in \cite{Despres04exp}. However, a uniform bound on the weaker norm (notice the weights $V_{KL}^n$)
$$
Q_h(u_0,t):=\sum_{K_n \in \mathcal{M}_N}\sum_{L\in\partial K_n^-}\delta
t |V_{KL}^n| |u^n_L-u^n_K|
$$
would be enough to prove the error estimate \refe{thmL1_1} by use of
the techniques of Kuznetsov \cite{Kuznetsov76} (and, therefore,
regardless whether the conservation law under consideration is
linear, the remark here is true). This uniform bound on $Q_h(u_0,t)$ is observed in numerical
experiments and is probably true. However, in this paper, we do
not prove this claim and it remains an open problem, even in the linear setting. 
\end{remark}

Interpolating between the $h^{1/2}$-error estimate in the
$L^\infty_tL^2_x$-norm of~\cite{Despres04nl} and the estimate of
Theorem~\ref{thmL1}, we obtain, under a strict CFL condition, a
$h^{1/2}$-error estimate in the $L^\infty_tL^p_x$-norm for data in
$W^{1,p}(\R^d)$, $1\leq p \leq 2$.

The analogue result for initial data in $W^{1,p}(\R^d)$, $2\leq p \leq
\infty$, requires different techniques and will be addressed in a subsequent paper.

\subsection{Sketch of the proof of Theorem~\ref{thmL1}}\label{skproof}

We only sketch the proof of the first part of the Theorem; besides we
do not exhibit the dependency on $t$ and $\xi$ of the error estimates (the
constants $C$ may depend on $t$, $T$ or $\xi$). 

To prove Theorem~\ref{thmL1}, we  use the fact that
the difference between the exact and the approximate solutions satisfies
the weak formulation (Lemma~\ref{formfaibleuh}):
\begin{eqnarray}
\int_{\R^d\times \R_+} (u-u_h)(\phi_t +V \cdot \nabla \phi) +\int_\R
(u-u_h)(\cdot,0) \phi (x,0)dx&=&(\mu_h+\nu_h)(\phi),\label{weakL2prep}
\end{eqnarray}
for every $\phi \in C^1_c(\R^d\times \R_+)$, where the error term $(\mu_h+\nu_h)(\phi)$ 
encompasses the
consistency error of the Finite Volume approximation 
and depends on the approximate solution.

We then notice that, by linearity of the equations, and subsequent principle of
superposition, we may suppose that the initial data $u_0\in
BV(\R^d)$ is the characteristic function $u_0=\mathbf{1}_A$ of a set $A$ with finite perimeter. 
The norm
$\|u_0\|_{TV}$ is equal to the perimeter of $A$. In fact, we show that we can also reduce the study to the case where $A$ does not
contain very thin parts so that, in particular, the volume of the
$h^{1/2}$-neighborhood of $\partial A$ 
\begin{eqnarray*}
A_0 &:=& \left\{x\in \R^d\::\: d(x,\partial A)\leq h^{1/2}\right\}.
\end{eqnarray*}
is bounded by $C |\partial A| h^{1/2}$.  At this stage, we build a function
$\phi_0 \in C^1(\R^d,[-1,1])$ such that $\phi_0 \equiv 1$ on $A\setminus
A_0$, $\phi_0 \equiv -1$ on  $A^c\setminus A_0$ and $\|\nabla
\phi_0\|_\infty \leq C {h^{-1/2}}$. Then we set $\phi(x,s)=\phi_0(X(x,s))\mathbf{1}_{[0,t]}(s)$ so that
$\phi_t+V\cdot \phi=0$ on $\R^d\times [0,t]$.
After a regularization process, the weak formulation yields   
\begin{eqnarray*}
  \|(u_h-u)(\cdot,t)\|_{L^1}&\leq & C|\partial A| h^{1/2} + |\mu_h(\phi)|+|\nu_h(\phi)|.  
\end{eqnarray*}
The term $\mu_h(\phi)$ is
\begin{eqnarray*}
\mu_h(\phi)&=&\sum_{K\in {\cal M}_N} |K| |u_K^{n+1}-u_K^n|
\left(\bbarint{K_n}{} \phi -\bbarint{K}{}\phi(\cdot,(n+1)\delta t)\right).
\end{eqnarray*}
The purpose of the reduction to the case $u_0=\mathbf{1}_A$ is an accurate estimation of the term $\mu(\phi)$. Indeed, by definition of $\phi$, we have: for $0\leq s \leq t$,  $\phi(x,s)= 1$ if $x \in X(A\setminus
A_0,s)$ and $\phi(x,s)= -1$ if $x \in X(A^c\setminus A_0,s)$. More
precisely, there exists $C>0$ such that, if 
$$d(K_n,\cup_{0\leq s \leq t}(X(\partial A,s) \times \{s\})  \geq C
(h^{1/2}+h+\delta t),$$ then $\phi$ is constant on the cell $K_n$. We denote by
$\mathcal{M}_N^1$ the set of cells which do not satisfy this
property; we have
\begin{eqnarray*}
\mu_h(\phi)&=&\sum_{K\in {\cal M}_N^1} |K| |u_K^{n+1}-u_K^n|
\left(\bbarint{K_n}{} \phi -\bbarint{K}{}\phi(\cdot,(n+1)\delta t)\right).
\end{eqnarray*}
and the $(d+1)$-dimensional volume of the cells of ${\cal M}_N^1$ is
bounded by $C |\partial A| h^{1/2}$ (Lemma \ref{idee1}). We use the Cauchy-Schwarz inequality and the estimate
$\|\partial_t \phi\|_\infty \leq {C}{h^{-1/2}}$ to get 
\begin{eqnarray*}
|\mu_h(\phi)|&\leq C E(u_0,t)^{1/2} |\partial A|^{1/2} h^{1/4}.    
\end{eqnarray*}
A similar technique yields the same bound on $\nu_h(\phi)$.

Thanks to the Energy estimate (Lemma~\ref{lemdeuxcasL2}) $E_h(u_0,t)
\leq C \mathcal{E}_h(u_0,t)$, we get 
\begin{eqnarray*}
 \|(u_h-u)(\cdot,t)\|_{L^1}&\leq & C\left(|\partial A| h^{1/2} \mathcal{E}_h(u_0,t)^{1/2} +|\partial A|^{1/2} h^{1/4}\right).  
\end{eqnarray*}
Finally, since the $L^2$-norm of the exact solution is
conserved, we get (Lemma \ref{idee2}): 
\begin{eqnarray*}
\mathcal{E}_h(u_0,t)&\leq &2 \|(u_h-u)(\cdot,t)\|_{L^1}.
\end{eqnarray*}
Together with the preceding inequality, this yields successively 
$\mathcal{E}_h(u_0,t) \leq C |\partial A| h^{1/2}$ and
$\|(u_h-u)(\cdot,t)\|_{L^1} \leq C |\partial A| h^{1/2}$.


\section{Classical results}\label{s2}

By continuity in $BV$ of the $L^2$-projection on the space of
functions which are constant with respect to the mesh ${\cal T}$, we
have:
\begin{lemma}\label{lemcihbv} There exists a constant $c\geq 0$ only depending on $d$ such that for every $u_0\in BV(\R^d)$,
\begin{eqnarray*}
\|u_h(\cdot,0)\|_{TV}= \frac{1}{2}\sum_{K\in \mathcal{M} ,\, L\in
  \partial K}
\!\!\!\!|K\!\!\shortmid\!L|\, |u_L^0-u_K^0|&\leq &\ds {c}{\alpha^{-3}} \, \|u_0\|_{TV},\\
\|u_h(\cdot,0)-u_0\|_{L^1}&\leq& \ds {c}{\alpha^{-1}}\,\|u_0\|_{TV} h.
\end{eqnarray*}
\end{lemma}

{\bf Proof of Lemma~\ref{lemcihbv}:} 
Using the weak density of $C^{1}_c(\R^d)$ in $BV(\R^d)$, we may
suppose that $u_0 \in C^{1}_c(\R^d)$. Let $K \in \mathcal{T}$, $L\in \partial K$. Since $|x-y|\leq 2h$ for
every $(x,y) \in K\times L$, we have
\begin{eqnarray*}
|u^0_K-u^0_L|&\leq &\cfrac{1}{|K|\,|L|} \int_{K}\int_L
|u_0(x)-u_0(y)| dxdy\\
&\leq & \cfrac{2h}{|K|\,|L|} \int_{K}\int_L \int_0^1 |\nabla
u_0((1-r)x+ry)| dr dx dy.
\end{eqnarray*} 
Now we perform the change of variables: $(x,y,r)\mapsto
(w=x-y,z=(1-t)x+ty,r=r)$ (of Jacobian determinant equal to $1$). We have
\begin{eqnarray*}
|u^0_K-u^0_L|&\leq & \cfrac{2h}{|K|\,|L|} 
\int_{B(x_K,2h)} |\nabla u_0(z)| \left(\int_0^1 \int_{\R^d} g(w,z,r) dw dr \right)dz,
\end{eqnarray*}
where $x_K\in \R^d$ is the centroid of $K$, and $g$ is defined by
$g(w,z,r)=1$ if $z+rw \in K$ and $z-(1-r)w \in L$, and $g(w,z,r)=0$ otherwise. For $(z,r)\in B(x_K,2h) \times
[0,1]$, we have $\int_{\R^d} g(w,z,r) dw \leq 2^d |K|$ if $r\geq 1/2$
and $\int_{\R^d} g(w,z,r) dw \leq 2^d |L|$ if $r < 1/2$. Finally, we obtain 
\begin{eqnarray*}
|u^0_K-u^0_L|&\leq & \cfrac{2^{d+1} \max (|K|,|L|)h}{|K|\,|L|}  \int_{B(x_K,2h)} |\nabla u_0(z)|dz.
\end{eqnarray*} 
We deduce from the inequalities~\eqref{eqgentilmesh} that for every
$K\in \mathcal{T}$,
\begin{eqnarray*}
\sum_{L\in \partial K}  |K\!\!\shortmid\!L||u^0_K-u^0_L|&\leq & 2^{d+1} \alpha^{-2} \int_{B(x_K,2h)}|\nabla u_0(z)|dz.
\end{eqnarray*}
Summing on $K\in \mathcal{T}$, we get 
\begin{eqnarray*}
\sum_{K\in \mathcal{T}}\sum_{L\in \partial K} |K\!\!\shortmid\!L|
|u^0_K-u^0_L| &\leq & 2^{d+1} \alpha^{-2}  \int_{\R^d} |\nabla u_0(z)|
M_{2h}(z) dz,
\end{eqnarray*}
where $M_{2h}(z)$ is the cardinal of the set $\{K\,:\, d(K,z)\leq
2h\}$. By~\eqref{eqgentilmesh}, we have 
\begin{eqnarray*}
M_{2h}(z) \alpha h^d &\leq  \sum_{K\,:\, d(K,z)\leq 2h} |K| &\leq |B(z,3h)|.
\end{eqnarray*} 
 Thus $ M_{2h}(z) \leq c \alpha^{-1}$ and we get the first part of
 the Lemma.

\medskip
Let $K\in \mathcal{T}$. Similarly, we have
\begin{eqnarray*}
\int_K |u_h(\cdot,0)-u_0| 
&\leq &  \cfrac{1}{|K|} \int_{K\times K} |u_0(x)-u_0(y)| dx dy\\
&\leq & {2^{d}h}\int_{B(x_K,h)} |\nabla u_0(z)|dz.
\end{eqnarray*} 
Summing on $K\in \mathcal{T}$ and using the fact that the cardinal of
the set $\{K \,:\,d(K,z)\leq h\}$ is bounded by ${c}{\alpha}^{-1}$, we get the
second estimate. 
\qed

The divergence free assumption~\refe{divergence} leads to the following identity
\begin{lemma}\label{conservation} Under the hypothesis on the
  divergence of $V$~\refe{divergence}, one has
\begin{equation*}
\sum_{L\in \partial K_n^-} V_{KL}^n= \sum_{L\,:\, K\in \partial L_n^-}
V_{LK}^n \qquad  \forall  \, K_n \in {\cal M}.
\end{equation*}
\end{lemma}

\subsection*{Monotony}\label{secmonotony}

Under a CFL condition the Finite Volume scheme is order-preserving:

\begin{proposition}\label{propmon} Under condition~\refe{CFL}, the
  linear application ${\cal L}:(u^n_K)\mapsto (u^{n+1}_K)$ defined by
  \refe{eqscheme} is order-preserving and stable for the $L^\infty$-norm.
\end{proposition}

{\bf Proof of Proposition~\ref{propmon}:} Eq.~\refe{eqscheme} gives
\begin{equation}\label{eqcc}
u^{n+1}_K=\left(1+\ds\sum_{L\in\partial K_n^-}\frac{V_{KL}^n\delta t}{|K|}\right)u^n_K-\ds\sum_{L\in\partial K^-_n}\frac{V_{KL}^n\delta
t}{|K|}u^n_L
\end{equation}
i.e., under \refe{CFL}, $u^{n+1}_K$ is a convex combination of
$u^n_K,(u^n_L)_{L\in\partial K_n^-}$. The stability result follows
from $\mathcal{L}(1)=1$.\qed

\subsection*{$L^1$-stability}
From Lemma~\ref{conservation} it is not difficult to see that the quantity $\sum_{K\in \cal T} |K| u_K^n$ is conserved. By the Crandall-Tartar Lemma~\cite{CrandallTartar80}, this fact and Proposition~\ref{propmon} imply 
\begin{proposition}
Under condition~\refe{CFL}, the scheme $\mathcal{L}:(u^n_K)\mapsto (u^{n+1}_K)$ is stable for the $L^1$-norm:
\begin{eqnarray*}
\sum_{K\in \cal T} |K| |u_K^{n+1}| &\leq &\sum_{K\in \cal T} |K| |u_K^{n}|,\qquad \forall n\geq 0.
\end{eqnarray*}
\label{l1stability}
\end{proposition}



\section{Weak Formulation}
\label{s3}
In the remainder of the paper, we assume that $u_0\in BV(\R^d)$, that the function $u\in C(\R_+, L^1(\R^d))$ is the exact solution to~\refe{PB}
and $u_h$ is the numerical approximation given by the
scheme~\eqref{eqscheme0}-\eqref{eqscheme}-\eqref{equh}. We also assume
that~\eqref{eqgentilmesh}, the CFL conditions
\eqref{CFL}-\eqref{CFL2} and the condition~\eqref{c1} are satisfied. 

We intend to prove that $u_h$ satisfies~\eqref{formfaible} up to an error term. 
\begin{lemma}
\label{lemformfaibleuh}
For every $\phi\in{\cal C}^\infty_c(\R^d \times[0,+\infty))$,
\begin{eqnarray}
\label{formfaibleuh}
 \ds\int_{\R_+}\int_{\R^d} u_h(\phi_t+V\cdot\nabla\phi)+\int_{\R^d} u_h(\cdot,0) \phi(\cdot,0) =
\mu_h(\phi)+\nu_h(\phi),
\end{eqnarray}
\begin{eqnarray}
\label{eqdefmu}
\mu_h(\phi)&:=&\sum_{K_n\in{\cal M}}|K|(u^{n+1}_K-u^n_K)(\langle\phi\rangle^n_K-\langle\phi\rangle_K((n+1)\delta t)),\\
\label{eqdefnu}
\nu_h(\phi)&:=&\sum_{K_n\in{\cal M}}\sum_{L\in\partial K_n^-}\delta
t(u^n_L-u^n_K)(V_{KL}^n\langle\phi\rangle^n_K- |K\!\!\shortmid\!L| \langle V
\cdot \mathbf{n} \,\phi\rangle^n_{KL}),
\end{eqnarray}
where 
\begin{eqnarray*}
\langle\phi\rangle^n_K:=\ds\bbarint{K_n}{} \phi\,,&\qquad \langle\phi\rangle_K(t):=\ds\bbarint{K}{} \phi(\cdot,t)\,,
&\qquad\langle V\cdot \mathbf{n}\,\phi\rangle^n_{KL}:=\ds \bbarint{K\!\shortmid L_n}{} \phi V \cdot \mathbf{n}_{KL}.
\end{eqnarray*}
\end{lemma}

{\bf Proof of Lemma~\ref{lemformfaibleuh}:} 
We develop the first term of the inequality and perform a discrete
integration by parts: 
\begin{multline}
\ds\int_0^\infty\int_{\R^d} u_h(\phi_t+V\cdot\nabla\phi)dx dt
=\ds\sum_{K_n\in{\cal M}}|K|\delta t\, u^n_K \bbarint{K_n}{}(\phi_t+V\cdot\nabla\phi)dxdt\\
=\ds\sum_{K_n\in{\cal M}}u^n_K\left(|K|\Big(\langle\phi\rangle_K((n+1)\delta
t)-\langle\phi\rangle_K(n\delta t) \Big)+\sum_{L\in\partial
  K}\delta t|K\!\!\shortmid\!L| \langle V\cdot \mathbf{n} \, \phi\rangle _{KL}^n\right)\\
=\ds \sum_{K_n \in{\cal M}} \left(|K|(u_K^n-u_K^{n+1})
\langle \phi\rangle_K((n+1)\delta t) + \sum_{L\in\partial K}\delta
t |K\!\!\shortmid\!L| \langle V\cdot \mathbf{n} \, \phi\rangle_{KL}^n u_K^n\right) \\
 - \sum_{K_n\in{\cal M}}
  |K|\delta t u_K^0 \langle\phi\rangle_K(0). 
\label{nimp0}
\end{multline}
The last sum is equal to $\int_{\R^d} u_h(\cdot,0)\phi(\cdot,0)$. Rearranging the terms in the second sum
and using the identity $\langle V \cdot\mathbf{n}\,\phi\rangle_{KL}^n+\langle V \cdot\mathbf{n}\,\phi\rangle_{LK}^n=0$, we
compute: for every $n\geq 0$
\begin{eqnarray*}
\sum_{K\in {\cal T}} \sum_{L \in \partial K} |K\!\!\shortmid\!L| \langle V \cdot\mathbf{n}\,\phi\rangle_{KL}^n u_K^n=
\sum_{K\in {\cal T}} \sum_{L \in \partial K_n^-} |K\!\!\shortmid\!L| \langle V \cdot\mathbf{n}\,\phi\rangle_{KL}^n(u_K^n- u_L^n).
\end{eqnarray*}
Plugging this identity and the definition of the scheme~\eqref{eqscheme}
in~\eqref{nimp0}, we get~\eqref{formfaibleuh}. \qed


We now assume that  $u_0\in BV(\R^d)$. Before building the test
function $\phi$, we show that we may reduce the study to more simple
initial data. 

Let ${\cal T}'$ be the cartesian mesh
$$\left\{t^{1/2}h^{1/2}\prod_{k=1}^d[p_k,p_k+1]\::\: (p_1,\cdots,p_d)\in
\mathbf{Z}^d\right\}.$$ 
The diameter of the cells is $h'=d^{1/2}t^{1/2}h^{1/2}$; this mesh is uniform: conditions~\eqref{eqgentilmesh} are satisfied for $\alpha'=\min(d^{-d/2}, d^{(d-3)/2}/2)$. In particular, by
  Lemma~\ref{lemcihbv}, we have  
\begin{eqnarray*}
\|v_0\|_{TV}\leq  C\|u_0\|_{TV} ,&\qquad& \|v_0-u_0\|_{L^1}\leq
 C\| u_0\|_{TV} t^{1/2} h^{1/2}, 
\end{eqnarray*}
where $v_0$ is the projection defined by  $\ds v_0 := \sum_{K'\in
  {\cal T}'} \left(\bbarint{\;K'}{} u_0\right)  \mathbf{1}_{K'}$.
These estimates, the $L^1$-stability of the scheme and of the
  evolution equation~\eqref{PB} imply that it is sufficient to prove Theorem~\ref{thmL1} for the initial data $v_0$. 
\begin{lemma}\cite{Federer69,Brenier84}
\label{lemdecompositionBV}
Let $v_0 \in BV(\R^d)$. For every $\eta \in \R$, define 
\begin{eqnarray*}
\chi_{v_0} (x,\eta) &:=& \left\{
\begin{array}{rl}
1 & \mbox{ if }\quad 0<\eta<v_0(x),\\
-1 & \mbox{if }\quad  v_0(x)<\eta <0,\\
0 & \mbox{ in the other cases}.
\end{array}
\right.
\end{eqnarray*}
We have $v_0=\ds\int_\R\chi_{v_0} (\cdot,\eta)  d\eta$, almost everywhere. Moreover,
\begin{eqnarray*}
\|v_0\|_{L^1}&=& \int_\R \|\chi_{v_0} (\cdot,\eta) \|_{L^1} d\eta,\\
\|v_0\|_{TV}&=& \int_\R \|\chi_{v_0} (\cdot,\eta) \|_{TV} d\eta,
\end{eqnarray*}
where the second identity is a consequence of the co-area formula for $BV$ functions.
\end{lemma}
Applying this decomposition, we deduce from the
linearity of the scheme and of the initial problem, that it is sufficient to
prove the Theorem for any initial data $u_0$ which is a finite sum of
characteristic functions $\mathbf{1}_{K'}$.  

Therefore, from now on, we suppose that 
$u_0=\mathbf{1}_A$ with $A=\cup_{1\leq j \leq l} K'_j$ and $K'_1,\cdots,K'_l$ are distinct elements of
${\cal T}'$. 


Let us point out that
$$
\|u_0\|_{L^1} = |A|=l t^{d/2}h^{d/2},
$$
and since, the boundary $\partial A$ is a finite union of distinct edges
denoted $K_1'\!\!\shortmid\!L_1', \cdots, K_m'\!\!\shortmid\!L_m'$, we deduce from the definition of the Total
Variation:
$$
\|u_0\|_{TV} =|\partial A| =\sum_{i=1}^m  |K_i'\!\shortmid\!L_i'|= m t^{(d-1)/2} h^{(d-1)/2}. 
$$

With the notations of Theorem~\ref{thsolPB}, we have $u=\mathbf{1}_B$ where
$B\subset \R^d\times \R_+$ is defined by 
\begin{eqnarray*}
B&:=&\ds\bigcup_{s\geq 0} Y(A,s)\times \{s\}.
\end{eqnarray*}

We are now going to build a test function $\phi$
for~\eqref{formfaible}~and~\eqref{formfaibleuh}. Let $\Gamma \in
C^\infty(\R,\R)$ such that $0\leq \Gamma \leq 1$, $\Gamma(R)=0$ for
$R\leq 1/3$ and $\Gamma(R)=1$ for $R\geq 2/3$. We define $\Phi_0 \in
\Lip(\R^d,\R)$ by 
\begin{eqnarray*}
\Phi_0(x)&:=&\left\{
\begin{array}{rcl}
\Gamma(t^{-1/2}h^{-1/2}d(x,\partial A)) &\mbox{ if }&x\in A,\\
-\Gamma(t^{-1/2}h^{-1/2}d(x,\partial A)) &\mbox{ if }&x\in A^c.
\end{array}\right.
\end{eqnarray*}
Since $d(\cdot,\partial A)$ is a $1$-Lipschitz continuous function, we have $\|\nabla \Phi_0\|_\infty
\leq t^{-1/2}h^{-1/2}\|\Gamma'\|_\infty$. Let us introduce a mollifier
$\rho \in C^\infty_c(\R^d,\R_+)$ such that $\int_{\R^d} \rho =1$ and $\supp \rho
\subset B(0,1/3)$, we set $\phi_0:=t^{-d/2}h^{-d/2}\rho(t^{-1/2}h^{-1/2}\,\cdot\,)\star
\Phi_0$. The function $\phi_0$ belongs to $C^1(\R^d)$ with the bound:
\begin{eqnarray*}
\|\nabla \phi_0\|_\infty &\leq & t^{-1/2}h^{-1/2}\|\Gamma'\|_\infty.
\end{eqnarray*}
Moreover, we have
\begin{eqnarray*}
x\in A \Longrightarrow \phi_0(x) \geq 0,&\quad \mbox{ and}\quad&x \in
A^{c}   \Longrightarrow \phi_0(x) \leq 0,
\end{eqnarray*}
and if we split $\R^d$ in three disjoint subsets:
$\R^d=\phi_0^{-1}(\{1\}) \cup \phi_0^{-1}((-1,1))\cup
\phi_0^{-1}(\{-1\})=: A_+\cup A_0\cup A_-$, we have $u_0\equiv 1$ on
$A_+$, $u_0 \equiv 0$ on $A_-$ and 
\begin{eqnarray*}
A_0&\subset & \left\{ x\;:\; d(x,\partial A)<
t^{1/2}h^{1/2}\right\} \;\subset \; \cup_{i=1}^m
B(x_i,{1}/{2}\,t^{1/2}h^{1/2}),
\end{eqnarray*}
where $x_i\in \R^d$ is the centroid of $ K_i'\!\!\shortmid\!L_i'$,
$1\leq i \leq m$. Consequently the Lebesgue measure of $A_0$ satisfies:
\begin{eqnarray}
\label{A0}
  |A_0|&\leq  C m t^{d/2}h^{d/2} &= C \|u_0\|_{TV} t^{1/2}h^{1/2} 
\end{eqnarray}
Finally, we define our test-function by $\phi(x,s):=\phi_0(X(x,s))
\mathbf{1}_{[0,t]}(s)$ for $(x,s)\in \R^d \times \R_+$, where $X$ is
the mapping introduced in Theorem~\ref{thsolPB}. Since $\div
V(\cdot,s) \equiv 0$ for every $s\geq 0$, we have
\begin{eqnarray*}
\phi_t + V\cdot \nabla \phi &=\phi_t + \div (V \phi)&=0 \qquad \mbox{in } \R^d\times (0,t), 
\end{eqnarray*}
and from~the estimate on $\nabla \phi_0$ above and on $\nabla X$ and $\partial_t X$
\eqref{estimationX} in Theorem~\ref{thsolPB}, $\phi$ belongs to $C^1(\R^d\times[0,t])$ and satisfies
\begin{eqnarray}
\label{nablaphi}
\|\nabla \phi \|_{L^\infty(\R^d\times [0,T])}&\leq & CC_0 t^{-1/2} h^{-1/2},\\
\label{dtphi}
\|\partial_t \phi \|_{L^\infty(\R^d\times [0,T])} &\leq & CC_0 t^{-1/2} h^{-1/2}.
\end{eqnarray}
Moreover, since $u(x,s)=u_0(X(x,s))$, $\phi$ inherits the properties
of $\phi_0$: for $(x,s)\in \R^d\times\R_+$,
\begin{eqnarray}
\label{phi2} 
u(x,s)=1 \Longrightarrow \phi(x,s) \geq 0,\quad
\mbox{ and }\quad  u(x,s)=0 \Longrightarrow \;\phi(x,s) \leq
0. 
\end{eqnarray}

\begin{lemma}
\label{lemmollifying}
We have
\begin{eqnarray}
\label{estimeL1_1}
\int_{\R^d} (u-u_h)(\cdot,t)\phi(\cdot,t) -\int_{\R^d}
(u-u_h)(\cdot,0)\phi(\cdot,t)&=&\mu_h(\phi)+\nu_h(\phi).
\end{eqnarray}
\end{lemma}
{\bf Proof of Lemma~\ref{lemmollifying}:}
The test function $\phi$ is not compactly supported and not time differentiable
at the time $t$. So it can not be used directly in~\eqref{formfaible},~\eqref{formfaibleuh}. First, remark that, since $u_0$ is compactly supported,
there exists $R>0$ such that $\supp u_h(\cdot,s) ,\supp u(\cdot,s)
\subset B(0,R)$ for every $0\leq s \leq t+\delta t$, thus $\phi$ may
be replaced in~\eqref{formfaible},~\eqref{formfaibleuh} by the compactly supported function $\phi \cdot \chi$
where $\chi \in C^\infty_c(\R^d)$ is such that $\chi \equiv 1$ on $B(0,R)$.  

To overcome the non differentiability at $t$ we introduce a mollifying
sequence to approximate the function $\mathbf{1}_{[0,t]}$.  
Let $(\psi_q)_{q\geq1}$ be a sequence of $C^\infty_c(\R_+)$ functions satisfying $\psi_q \equiv 1$ on $[0,t]$,
$\psi_q\equiv 0$ on $[t+{\delta t}/{q},+\infty)$ and $0\geq \psi_q'\geq
  {-2q}/{\delta t}$. We  set $\phi_{q}(x,s) := \phi_0(X(x,s)) \psi_q(s)$ for every $(x,s) \in
  \R^d\times \R_+$ and $q\geq 1$. The weak formulation~\eqref{formfaibleuh} reads
\begin{eqnarray*}
\ds\int_{t}^{t+\delta t} \left(\int_{\R^d} u_h(x,s) \phi_0(X(x,s))
dx\right) \psi'_q(s) ds +\int_{\R^d} u_h(\cdot,0) \phi_0(\cdot,0) 
=\mu_h(\phi_{q}) +\nu_h(\phi_{q}),
\end{eqnarray*}
The sequence of Radon measures $(\psi'_q(s)ds)_q$ converges to
  $-\delta_{t}$, so the first term of the left hand side converges towards
$$ \ds -\int_{\R^d} u_h(\cdot,t) \phi(\cdot,t)$$ 
as $q$ tends to $+\infty$.

Clearly, for every $K_n\in
  \mathcal{M}$ and  $L\in \partial
K_n^-$, we have $(\langle \phi_q \rangle_K)_q^n
\rightarrow \langle \phi \rangle_K^n$,  $(\langle V\cdot \mathbf{n} \phi_q \rangle_{KL}^n)_q
\rightarrow \langle V\cdot \mathbf{n}  \phi \rangle_{KL}^n$ and $(\langle \phi_q
\rangle_{K}(n\delta t))_q
\rightarrow \langle \phi\rangle_K(n\delta t)$. 
Passing to the limit on $q$, we get
\begin{eqnarray*}
\ds-\int_{\R^d} u_h(\cdot,t) \phi(\cdot,t) +\int_{\R^d} u_h(\cdot,0) \phi_0 
=\mu_h(\phi) +\nu_h(\phi).
\end{eqnarray*}
Similarly, we obtain from~\eqref{formfaible} 
\begin{eqnarray*}
\ds-\int_{\R^d} u(\cdot,t)\phi(\cdot,t)+\int_{\R^d} u(\cdot,0)\phi_0&=&0. 
\end{eqnarray*}
Subtracting these equalities, we get the result. \qed

We relate the first term of the equality of Lemma~\ref{lemmollifying} to the $L^\infty_t L^1_x$-norm of the error:
\begin{lemma} We have
\label{lemclaim}
\begin{eqnarray}
\label{claim} 
\|(u_h-u)(\cdot,t)\|_{L^1} &\leq & \int_{\R^d}
(u-u_h)(\cdot,t)\phi(\cdot,t) + C \|u_0\|_{TV} t^{1/2}h^{1/2} .
\end{eqnarray}
\end{lemma}
{\bf Proof of Lemma~\ref{lemclaim}:}
 Recall that from the monotony of the scheme, we have
\begin{eqnarray}
\label{eqmonotony}
0\leq u_h(x,s)&\leq 1, \qquad\forall \, (x,s)\in \R^d\times \R_+.
\end{eqnarray}
Let us consider the disjoint decomposition $\R^d=X(A_+,t) \cup
X(A_0,t) \cup X(A_-,t)$ and notice
that~\eqref{phi2},~\eqref{eqmonotony} imply that $(u(x,t)-u_h(x,t))\phi(x,t) \geq 0$ for
every $x\in \R^d$. In particular
\begin{eqnarray*}
0&\leq &\int_{X(A_0,t)} (u-u_h)(\cdot,t)\phi(\cdot,t).
\end{eqnarray*}
Since $0\leq u,u_h\leq 1$, we deduce from~\eqref{A0} and the
conservation property~\eqref{conservemesure} that
\begin{eqnarray*}
\int_{X(A_0,t)} |(u_h-u)(\cdot,t)|\,\leq&  |X(A_0,t)|=|A_0| &\leq  C \|u_0\|_{TV} t^{1/2}h^{1/2}.
\end{eqnarray*}
Finally, using the identities $\phi(\cdot,t)\equiv 1$,
$u(\cdot,t)\equiv 1$ on $X(A_+,t)$ and $\phi(\cdot,t)\equiv -1$,
$u(\cdot,t) \equiv 0$ on $X(A_-,t)$ and
the bounds $0\leq u_h \leq 1$, we have
\begin{eqnarray*}
\int_{X(A_+\cup A_-,t)} |u_h(\cdot,t) -u(\cdot,t)|&=& \int_{X(A_+,t)}
(u-u_h)(\cdot,t)+ \int_{X(A_-,t)} (u_h-u)(\cdot,t)\\& =&\int_{X(A_+\cup A_-,t)} (u-u_h)(\cdot,t)\phi(\cdot,t).
\end{eqnarray*}
Summing up this equality and the two preceding estimates, we get~\eqref{claim}. \qed
Lemmas~\ref{lemmollifying} and~\ref{lemclaim} and the second
estimate of Lemma~\ref{lemcihbv} (and $h\leq c_1 t$) yield 
\begin{eqnarray}
\label{estimeL1_2}
\|(u_h-u)(\cdot,t)\|_{L^1} &\leq
&\mu_h(\phi)+\nu_h(\phi)  + C \|u_0\|_{TV}\,t^{1/2}h^{1/2}.
\end{eqnarray}
The first key idea to estimate $\mu_h(\phi)$ and $\nu_h(\phi)$ is to
remark that most of the terms in the sums of~\eqref{eqdefmu},~\eqref{eqdefnu} vanish.
\begin{lemma}
\label{idee1}
There exists $\mathcal{M}_N^1 \subset \mathcal{M}_N$ such that 
\begin{eqnarray*}
\mu_h(\phi)&=&\sum_{K_n\in{\cal M}_N^1}|K|(u^{n+1}_K-u^n_K)(\langle\phi\rangle^n_K-\langle\phi\rangle_K((n+1)\delta t)),\\
\nu_h(\phi)&=&\sum_{K_n\in{\cal M}_N^1}\sum_{L\in\partial K_n^-}\delta
t(u^n_L-u^n_K)(V_{KL}^n\langle\phi\rangle^n_K- |K\!\!\shortmid\!L| \langle
V \cdot\mathbf{n}\, \phi\rangle^n_{KL}).
\end{eqnarray*}
with the estimate
\begin{eqnarray}
\label{volume}
\sum_{K_n\in \mathcal{M}_N^1 }\delta t |K| &\leq &  CC_0^d \|u_0\|_{TV}t^{3/2}h^{1/2}.
\end{eqnarray}
\end{lemma}
{\bf Proof of Lemma~\ref{idee1}:}
Let us denote by $\mathcal{B}$ the open subset of $\R^d \times [0,t]$ of all the points 
$(x,s)$ such that $|\phi(x,s)|<1$. Let $K_n\in {\cal
  M}_N$ which does not intersect $\mathcal{B}$. In this case,
  $\phi$ is constant ($\phi \equiv 1$ or $\phi\equiv-1$) on
 ${K_n}$ and we have
\begin{eqnarray*}
\langle\phi \rangle_K^n=\langle \phi \rangle_K((n+1)\delta t) \quad\mbox{ and } \quad V_{KL}^n\langle
\phi \rangle_{KL}^n =|K\!\!\shortmid\!L|\langle  V \cdot\mathbf{n}\,\phi \rangle_{KL}^n ,\qquad
\forall\, L\in \partial K_n^-.
\end{eqnarray*}
Consequently, the  corresponding terms in the definitions of
$\mu_h$ and $\nu_h$ vanish and the equalities of the Lemma hold with
\begin{eqnarray*}
\mathcal{M}_N^1&:=&\left\{K_n \in \mathcal{M}_N\;:\; K_n \cap
\mathcal{B} \neq \emptyset\right\},
\end{eqnarray*}
We now estimate the $(d+1)$-dimensional
Lebesgue measure of the union of the cells $K_n \in {\cal M}_N^1$.
From the definition of $\phi$, for every $K_n\in {\cal M}_N^1$, there
exist $(x,s) \in K_n$ such that  $d(X(x,s),\partial A) \leq
t^{1/2}\,h^{1/2}$. Since $\partial A=\cup_{1\leq i \leq m}
K'_i\!\shortmid \!L'_i$, for every $K_n\in {\cal M}_N^1$, there
exist $(x_{K,n},s_{K_n}) \in K_n$ and $1\leq i \leq m$ such that
$d(X(x_{K_n},s_{K_n}),x_i)\leq 2 t^{1/2}\,h^{1/2}$, where $x_i$ is the
centroid of $K'_i\!\shortmid \!L'_i$. Thus using the estimates
of~\eqref{estimationY} on $\nabla Y$ and $\partial_t Y$, we obtain 
 \begin{eqnarray*}
\bigcup_{K_n\in \mathcal{M}_N^1}(x_{K_n},{s_{K_n}}) &\subset & \ds \bigcup_{i=1}^m \bigcup_{0\leq s\leq t} Y(B(x_i,2t^{1/2}h^{1/2}),s) \times \{s\}.\\
&\subset &\bigcup_{i=1}^m \bigcup_{n=0}^N B\left(Y(x_i,n\delta t),2C_0 t^{1/2}\,h^{1/2} +C_0\delta t) \right)  \times  [n\delta t,(n+1)\delta t],
\end{eqnarray*}
and 
 \begin{eqnarray*}
\bigcup_{K_n\in \mathcal{M}_N^1}{K_n} &\subset &\bigcup_{i=1}^m \bigcup_{n=0}^N B\left(Y(x_i,n\delta t),2C_0 t^{1/2}\,h^{1/2} +C_0\delta t +h) \right)  \times  [n\delta t,(n+1)\delta t].
\end{eqnarray*}
Finally, using $C_0\geq 1$, $\delta t \leq c_0 h$ and $h \leq c_1 t$, we get
\begin{eqnarray*}
\sum_{K_n\in \mathcal{M}_N^1 }\delta t |K| \,\leq & CC_0^d m t
(t^{d/2}h^{d/2}+ h^{d})\,\leq &  CC_0^d \|u_0\|_{TV} t^{3/2} h^{1/2}.~~~~~~\mbox{\qed}
\end{eqnarray*}

\section{Energy Estimates}
\label{s4}
The two next lemmas contain the second key idea of the paper. The
first one shows that the $L^1$-error at time $t$ controls ${\cal
  E}_h(u_0,t)$ (the jump of the $L^2$-energy of the approximate
solution). The second one relates ${\cal E}_h(u_0,t)$ to the
regularity of the approximate solution. 
\begin{lemma}
\label{idee2}
We have
\begin{eqnarray*}
\mathcal{E}_h(u_0,t) &\leq& 2
\,\|(u_h-u)(\cdot,t)\|_{L^1}.
\end{eqnarray*}
\end{lemma}
{\bf Proof of Lemma~\ref{idee2}:}
Using the conservation property of Corollary~\ref{conservenormes} with $f(v)=v^2$, we have $
\|u(\cdot,t)\|_{L^2}^2 =\|u_0\|_{L^2}^2$. Thus
\begin{eqnarray*}
  \mathcal{E}_h(u_0,t) &=&\left(\|u_h(\cdot,0)\|_{L^2}^2  -\|u_0\|_{L^2}^2\right)
 +\left( \|u(\cdot,t)\|_{L^2}^2 - \|u_h(\cdot,t)\|_{L^2}^2\right).
\end{eqnarray*}
Since $u_h(\cdot,0)$ is the $L^2$-projection of $u_0$ on the mesh
$\mathcal{T}$, the first term is non
positive. Finally, $0\leq u,u_h\leq 1$ yields $u^2(\cdot,t)-u_h^2(\cdot,t) \leq 2\,|u_h(\cdot,t)-u(\cdot,t)|$.\qed

In particular,~\eqref{estimeL1_2} implies 
\begin{eqnarray}
\label{estimeL1_3}
 \mathcal{E}_h(u_0,t)&\leq &
2 \mu_h(\phi)+2\nu_h(\phi) + C \|u_0\|_{TV}\,t^{1/2}h^{1/2}.
\end{eqnarray}

\begin{lemma}
\label{lemE}
We have the following identity
\begin{multline}
\label{idlemE}
\mathcal{E}_h(u_0,t)=\sum_{K_n \in \mathcal{M}_N}\sum_{L\in \partial K_n^-}\!\!
\left(1+\sum_{M \in \partial K_n^-} \cfrac{V_{KM}^n \delta t}{|K|}\right) |V_{KL}^n|\delta
  t \,\left(u_L^n -u_K^n\right)^2\\
 +\cfrac{1}{2} \sum_{K_n \in \mathcal{M}_N} \sum_{L,M\in \partial
   K_n^-}
 \cfrac{V_{KL}^nV_{KM}^n\delta
  t^2}{|K|} \left(u_M^n -u_L^n\right)^2.
\end{multline} 
\end{lemma}
\begin{remark}
We say that the identity of
Lemma~\ref{lemE} is an Energy estimate, because of the analogy with
the \textit{a priori} inequality 
\begin{eqnarray*}
\|u_{\nu}(\cdot,0)\|_{L^2}^2-\|u_{\nu}(\cdot,t)\|_{L^2}^2&\leq& 2\nu
\int_0^t \int_{\R^d} |\nabla u_\nu(x,s)|^2 dx ds,
\end{eqnarray*}     
where $u_\nu$ satisfies the advection-diffusion equation $\partial_t
u_\nu +\div (V u_\nu) -\nu \Delta u_\nu=0$. In our case the diffusion is due
to the scheme (and does not have the regularity and the isotropy of the
Laplace operator), and the corresponding parameter $\nu$ is of order $h$.
\end{remark}

{\bf Proof of Lemma~\ref{lemE}:} Let $n\geq 0$. Using
Lemma~\ref{conservation} to change the index of summation in the last sum of the right hand-side of the first equality, we have
\begin{eqnarray*}
\sum_{K \in{\cal T}} |K|(u_K^n)^2&=&\sum_{K \in {\cal T} }
\left(|K|+\sum_{L\in \partial K_n^-} V_{KL}^n \delta t\right)
(u_K^n)^2 -\sum_{K\in {\cal T}} \sum_{L\::\:K\in \partial L_n^-}
V_{LK}^n \delta t (u_K^n)^2\\
&=& \sum_{K \in {\cal T} }
\left\{ \left(|K|+\sum_{L\in \partial K_n^-} V_{KL}^n \delta t\right)
(u_K^n)^2  - \sum_{L\in \partial K_n^-}  V_{KL}^n \delta t
(u_L^n)^2\right\}\\
&=&\sum_{K \in {\cal T}} |K|\sum_{L\in \mathcal{V}_n^-(K)} a_{KL}^n (u_L^n)^2,
\end{eqnarray*}
where we have set for every $K_n\in {\cal M}$, $\mathcal{V}_n^-(K):=\partial K_n^-\cup\{K\}$, and for every  $L \in \partial K_n^-$:
\begin{eqnarray*}
\quad
 a_{KK}^n:=1+\ds\sum_{M\in \partial K_n^-} \frac{V_{KM}^n\delta t}{|K|} ,&\quad &a_{KL}^n
:=-\ds\frac{V_{KL}^n\delta t}{|K|}.
\end{eqnarray*} 
With this notations,~\eqref{eqscheme} reads 
\begin{eqnarray*}
u_K^{n+1}=\sum_{L\in
  \mathcal{V}_n^-(K)} a_{KL}^n u_K^n,
\end{eqnarray*}
for every $K_n \in \mathcal{M}$, and we have
\begin{multline*}
\sum_{K \in{\cal T}} |K|\left((u_K^n)^2-(u_K^{n+1})^2\right)=
\sum_{K\in{\cal T}} |K|\left\{
\sum_{L\in \mathcal{V}_n^-(K)} a_{KL}^n (u_L^n)^2 -\left(\sum_{L\in\mathcal{V}_n^-(K)} a_{KL}^n u_L^n\right)^2 
\right\}\\
\qquad=\sum_{K\in{\cal T}} |K|\left\{\sum_{L\in\mathcal{V}_n^-(K)}
a_{KL}^n (1-a_{KL}^n)(u_L^n)^2-\sum_{L,M\in\mathcal{V}_n^-(K),\,M\neq L}
a_{KL}^na_{KM}^nu_L^nu_M^n
\right\}.\\
\end{multline*}
Using the identity $\sum_{L\in\mathcal{V}_n^-(K)} a_{KL}^n=1$, we compute
\begin{eqnarray*}
\sum_{K \in{\cal T}} |K|((u_K^n)^2-(u_K^{n+1})^2)&=&
\sum_{K\in{\cal T}} |K| \left\{\sum_{L\in\mathcal{V}_n^-(K)}a_{KL}^n u_L^n \left(\sum_{M\in\mathcal{V}_n^-(K)}
a_{KM}^n (u_L^n-u_M^n)\right)\right\}\\
&=&\cfrac{1}{2}\sum_{K\in{\cal T}}|K|\sum_{L,M\in\mathcal{V}_n^-(K)} a_{KL}^n a_{KM}^n (u_L^n-u_M^n)^2.
\end{eqnarray*} 
Summing on $0\leq n\leq N$, we get the Lemma. \qed

\section{Proof of the first part of Theorem~\ref{thmL1}}
\label{s41}
In this section, we assume that $\xi >0$. We estimate 
successively $\mu_h(\phi)$ and $\nu_h(\phi)$. Applying the
Cauchy-Schwarz inequality to the first formula of
Lemma~\ref{idee1}, we obtain
\begin{eqnarray}
\nonumber
|\mu_h(\phi)|&\leq & E_h(u_0,t)^{1/2} \left(\sum_{K_n\in{\cal M}_N^1}
|K| \delta t^2 \left( \cfrac{\langle \phi \rangle^n_K-\langle \phi \rangle_K((n+1)\delta
  t)}{\delta t}\right)^2 \right)^{1/2}.
\end{eqnarray}
On the basis of \refe{dtphi} and the CFL condition~\refe{CFL2} $\delta t\leq c_0 h$, we have
$$
|\langle \phi \rangle^n_K-\langle \phi \rangle_K((n+1)\delta t)|\leq \|\partial_t\phi\|_{L^\infty}\delta t\leq CC_0 \,\ds{\delta t}\,{t^{-1/2}h^{-1/2}}\leq CC_0\ds{h^{1/2}}{t^{-1/2}}
  $$
and therefore, by \refe{volume}
\begin{eqnarray}
\nonumber
|\mu_h(\phi)|&\leq & CC_0^{1+d/2}  E_h(u_0,t)^{1/2} \|u_0\|_{TV}^{1/2} \, t^{1/4} h^{1/4}.\label{etdedeux}
\end{eqnarray}
We now write
\begin{multline}
\nu_h(\phi)=\sum_{K_n\in{\cal M}_N^1}\sum_{L\in\partial K_n^-}\delta
t(u^n_L-u^n_K)V_{KL}^n\left(\langle\phi\rangle^n_K- \langle
\phi\rangle^n_{KL}\right)\\
+\sum_{K_n\in{\cal M}_N^1}\sum_{L\in\partial K_n^-}\delta
t(u^n_L-u^n_K) |K\!\!\shortmid\!L| \left(\langle V\cdot \mathbf{n}
\rangle^n_{KL}\langle\phi\rangle^n_{KL} -\langle
V \cdot\mathbf{n}\,\phi\rangle^n_{KL}  \right),\\
\qquad\qquad=: \mathrm{I} +\mathrm{II}.\hfill~ \label{IIIetiv}
\end{multline}
From the Cauchy-Schwarz inequality, we have
\begin{eqnarray*}
 |\mathrm{I}|&\leq & E_h(u_0,t)^{1/2} \left(\sum_{K_n\in{\cal
     M}_N^1}\sum_{L\in\partial K_n^-} |V_{KL}^n| \delta t h^2
  \left(\cfrac{\langle\phi\rangle^n_K- \langle
\phi\rangle^n_{KL} }{h}\right)^2\right)^{1/2}.
\end{eqnarray*}
The conditions~\eqref{eqgentilmesh} implies
\begin{eqnarray*}
\sum_{L \in \partial K_n^-} |V_{KL}^n| h&\leq  \alpha^{-2}
\|V\|_\infty |K|&\leq C |K|,\qquad \forall
\, K_n \in \mathcal{M}.
\end{eqnarray*}
Consequently,~\eqref{nablaphi} and~\eqref{volume} lead to
\begin{eqnarray}
\label{etdetrois}
 |\mathrm{I}|&\leq &CC_0^{1+d/2} E_h(u_0,t)^{1/2} \|u_0\|_{TV}^{1/2}\, t^{1/4} h^{1/4}.
\end{eqnarray}
We now estimate the term $\mathrm{II}$; notice that it vanishes in case $V={\rm Cst}$.
Let $K_n \in \mathcal{M}_N^1$ and $L\in \partial K_n^{-}$, since
$\ds \bbarint{K\!\shortmid L_n}{} \langle V \cdot\mathbf{n}\rangle_{KL}^n - V\cdot
\mathbf{n}_{KL} =0$, we have
\begin{multline*}
\langle V\cdot \mathbf{n}
\rangle^n_{KL}\langle\phi\rangle^n_{KL} -\langle
V \cdot\mathbf{n}\,\phi\rangle^n_{KL}= \bbarint{K\!\shortmid L_n}{} \phi(x,r)
\left( \langle V \cdot\mathbf{n}\rangle_{KL}^n - V(x,r)\cdot
\mathbf{n}_{KL} \right) dxdr\\
= \bbarint{K\!\shortmid L_n}{} \left(\phi(x,r) -\langle \phi\rangle_{KL}^n \right)
\left( \langle V \cdot\mathbf{n}\rangle_{KL}^n - V(x,r)\cdot
\mathbf{n}_{KL} \right) dxdr.
\end{multline*}
Using~\eqref{nablaphi}-\eqref{dtphi} and $|\langle V \cdot\mathbf{n}\rangle_{KL}^n - V(x,r)\cdot
\mathbf{n}_{KL}|\leq Ch$ for $(x,r) \in K\!\!\shortmid\!L_n$, we get
\begin{eqnarray*}
\left|\langle V\cdot \mathbf{n}
\rangle^n_{KL}\langle\phi\rangle^n_{KL} -\langle
V \cdot\mathbf{n}\,\phi\rangle^n_{KL}\right|& \leq & CC_0
\,t^{-1/2}h^{3/2}. 
\end{eqnarray*} 
Thus, we have
\begin{eqnarray*}
|\mathrm{II}|&\leq & CC_0 \left(\sum_{K_n\in{\cal M}_N^1}\left(\sum_{L\in\partial
 K_n^-} |K\!\!\shortmid\!L|h\right) \delta t \right)t^{-1/2} h^{1/2}. 
\end{eqnarray*}
As above, for every $K_n \in \mathcal{M}$ we have $\sum_{L\in\partial
 K_n^-} |K\!\!\shortmid\!L|h \leq C |K|$ and from \eqref{volume}:
\begin{eqnarray}
\label{etdequatre}
|\mathrm{II}| &\leq & C C_0^{d+1} \|u_0\|_{TV} th.
\end{eqnarray}
Summing up~\eqref{etdedeux}-\eqref{etdetrois}-\eqref{etdequatre}, we obtain
\begin{eqnarray}
\label{munu}
|\mu_h(\phi)|+|\nu_h(\phi)| \leq  C \left( C_0^{d/2+1} E_h(u_0,t)^{1/2}
\|u_0\|_{TV}^{1/2}\, t^{1/4} h^{1/4} +  C_0^{d+1}\|u_0\|_{TV} t h \right),
\end{eqnarray}
and~\eqref{estimeL1_3} yields
\begin{eqnarray}
\label{estimeL1_4}
  \mathcal{E}_h(u_0,t) \leq  C  \left( C_0^{d/2+1} E_h(u_0,t)^{1/2}
\|u_0\|_{TV}^{1/2}\, t^{1/4} h^{1/4} + \|u_0\|_{TV}(t^{1/2} h^{1/2}+  C_0^{d+1}th) \right).
\end{eqnarray}
The following Lemma will allow us to bound $E_h(u_0,t)$ by $\mathcal{E}_h(u_0,t)$.
\begin{lemma}
\label{lemdeuxcasL2}
\begin{eqnarray*}
E_h(u_0,t)&\leq &\ds{C}{\xi}^{-1}\, \mathcal{E}_h(u_0,t).
\end{eqnarray*}  
\end{lemma}

{\bf Proof of Lemma~\ref{lemdeuxcasL2}:}
From~\eqref{eqscheme} and the Cauchy-Schwarz inequality, we have for
every $K_n\in{\cal M}$,
\begin{eqnarray*}
|K| |u_{K}^{n+1}-u_K^n|^2&\leq &\left(\sum_{L\in \partial K_n^-}
 \ds\frac{|V_{KL}^n|\delta t}{|K|}\right) \sum_{L\in \partial K_n^-}
 |V_{KL}^n|\delta t |u_L^{n}-u_K^n|^2. 
\end{eqnarray*}
Due to the CFL condition~\refe{CFL}, we have $\sum_{L\in \partial K_n^-}
 {|V_{KL}^n|\delta t}\leq|K|$; consequently, summing the preceding inequality on $K_n \in {\cal
 M}_N$, we obtain
\begin{eqnarray*}
E_h(u_0,t)&\leq & 2 \sum_{K_n \in \mathcal{M}_N}\sum_{L \in
  \partial K_n^-} |V_{KL}^n|\delta t |u_L^n-u_K^n|^2.
\end{eqnarray*}
Using again the CFL condition,  for every
$K_n\in{\cal M}$ we have $(1+\sum_{M\in \partial K_n^-} V_{KM}\delta
t/|K|)\geq \xi$. 
Thus $\xi E_h(u_0,t)/2$ is bounded by the first sum of the
  right hand side of~\eqref{idlemE}.\qed

Plugging this inequality in~\eqref{estimeL1_4}, we get
\begin{eqnarray*}
  \mathcal{E}_h(u_0,t) \leq  C  \left( C_0^{d/2+1}
  \mathcal{E}_h(u_0,t)^{1/2} \xi^{-1/2}
\|u_0\|_{TV}^{1/2}\, t^{1/4} h^{1/4} + \|u_0\|_{TV}(t^{1/2} h^{1/2}+  C_0^{d+1}th) \right).
\end{eqnarray*}
We then use  the Young
inequality $ab\leq \frac{1}{2}a^2+ \frac{1}{2}b^2$ with 
$$
a=C^{-1/2}\mathcal{E}_h(u_0,t)^{1/2},\quad b=C^{1/2} C_0^{d/2+1} {\xi^{-1/2}}\|u_0\|^{1/2}_{TV} t^{1/4} h^{1/4}
$$ 
to derive the estimate
\begin{eqnarray*}
\mathcal{E}_h(u_0,t)&\leq& CC_0^{d+2}\|u_0\|_{TV}  ({\xi}^{-1} 
\,t^{1/2} h^{1/2} + th).
\end{eqnarray*}
Finally, from Lemma~\ref{lemdeuxcasL2} and the inequalities~\eqref{munu},~\eqref{estimeL1_2}, we compute successively
\begin{eqnarray*}
E_h(u_0,t) & \leq &  CC_0^{d+2}{\xi}^{-2} \|u_0\|_{TV}  (
\,t^{1/2} h^{1/2} + {\xi}\,th),\\
|\mu_h(\phi)|+|\nu_h(\phi)| & \leq &  CC_0^{d+2}\xi^{-1}\|u_0\|_{TV} (t^{1/2}h^{1/2} + \xi^{1/2}\,th),\\
\|(u_h-u)(\cdot,t)\|_{L^1} &\leq &
CC_0^{d+2}\xi^{-1}\|u_0\|_{TV} (t^{1/2}h^{1/2} + \xi^{1/2}\,th).~~~~~~~~~~~~\mbox{\qed}
\end{eqnarray*}

\section{Proof of the Second part of Theorem~\ref{thmL1}}
\label{s42}
We suppose that $V$ does not depend on the time variable. If
$$\sup_{K\in \mathcal{T}} |K|^{-1}\sum_{L\in \partial K_0^-} |V_{KL}^0|\delta t
< 2/3,$$ then the proof of the first case allows us to conclude; therefore we can suppose that we have $ |K^\star|\leq 2\sum_{L\in \partial
  K^{\star,-}_0}|V_{K^\star L}^0|\delta t$  for some $K^\star \in {\cal
  T}$. This cell $K^\star$ satisfying~\eqref{eqgentilmesh}, we deduce that the following inverse CFL condition holds:
\begin{eqnarray}
\label{invCFL2}
h&\leq & 2\alpha^{-2}\|V\|_\infty\delta t . 
\end{eqnarray}
We now bound $\mu_h(\phi)$ and $\nu_h(\phi)$. From Lemma~\ref{idee1}
and the estimate~\eqref{dtphi} on $\partial_t \phi$, we have
\begin{eqnarray*}
 |\mu_h(\phi)|&\leq &C \delta t (t^{-1/2}h^{-1/2}) \sum_{K_n\in \mathcal{M}_N^1} |K| |u_K^{n+1}-u_{K}^n|. 
\end{eqnarray*}
Since $V$ does not depend on $t$, we have $V_K^n=V_K^0$, for every
$K_n\in {\cal M}$. Consequently the discrete time derivative $v_K^n
:=\delta t^{-1}(u_K^{n+1}-u_K^n)$ satisfies the scheme~\eqref{eqscheme}. The $L^1$-stability of the scheme and
Lemma~\ref{lemcihbv} imply 
\begin{eqnarray*}
\sum_{K_n\in \mathcal{M}_N}|K|\left|u_K^{n+1}-u_{K}^n\right|&\leq &\ds\frac{t}{\delta t}
  \sum_{K\in \mathcal{T}} |K| |u_K^1-u_K^0|
\, \leq  t \sum_{K\in \mathcal{T}} |V_{KL}^0| |u_K^0-u_L^0| \\
&\leq &  t \|V\|_\infty \sum_{K\in \mathcal{T}}  |K\!\!\shortmid\!L|  |u_K^0-u_L^0|
  \,\leq  C \|u_0\|_{TV} t.
\end{eqnarray*}
Thus (using the CFL condition $\delta t\leq c_0h$), we have
\begin{eqnarray}
\label{un}
 |\mu_h(\phi)|&\leq &C \|u_0\|_{TV} t^{1/2} h^{1/2}.
\end{eqnarray}
As in the proof of the first case (see \eqref{IIIetiv}), we write $\nu_h(\phi)=\mathrm{I} +
\mathrm{II}$.
Recall that
\begin{eqnarray*}
\mathrm{I}=\sum_{K_n\in{\cal M}_N^1}\sum_{L\in\partial K_n^-}\delta
t(u^n_L-u^n_K)V_{KL}^n\left(\langle\phi\rangle^n_K- \langle
\phi\rangle^n_{KL}\right),
\end{eqnarray*}
and that, from the estimate~\eqref{nablaphi} on $\nabla \phi$, we have
 \begin{eqnarray*}
|\mathrm{I}|\leq CC_0\, h^{1/2}t^{-1/2}  \sum_{K_n\in{\cal M}_N^1}\sum_{L\in\partial K_n^-}\delta
t |V_{KL}^n ||u^n_L-u^n_K|.
\end{eqnarray*}
Using the definition of the scheme~\eqref{eqscheme}, for $K_n \in \mathcal{M}_N^1$ and $L\in
\partial K_n^-$, we have
\begin{eqnarray*}
u_L^n-u_K^n&=& \cfrac{|K|}{\sum_{M\in \partial
    K_n^-}V_{KM}^n\delta t}\left(u_K^{n+1} -u_K^n\right)+\cfrac{1}{\sum_{M\in \partial
    K_n^-}V_{KM}^n\delta t}\sum_{M\in \partial K_n^-}
V_{KM}^n \delta t(u_L^n-u_M^n).
\end{eqnarray*}
Therefore, 
\begin{multline*}
|\mathrm{I}|\leq  CC_0 \, h^{1/2}t^{-1/2} \Bigg(\sum_{K_n\in{\cal M}_N^1} |K||u_K^{n+1}-u_K^n| \\+\sum_{K_n\in{\cal M}_N^1} \cfrac{1}{\sum_{M\in \partial K_n^-}
|V_{KM}^n| \delta t}\sum_{L,M\in \partial K_n^-} V_{KL}^nV_{KM}^n \delta t^2
|u_M^n-u_L^n|\Bigg)\\
~\qquad~~~=CC_0\, h^{1/2}t^{-1/2}\, (\mathrm{I}_1 + \mathrm{I}_2).\hfill~
\end{multline*}
Using again the $L^1$-stability of the scheme, we have
\begin{eqnarray*}
|\mathrm{I}_1|&\leq & C\|u_0\|_{TV} \,t . 
\end{eqnarray*}
For $\mathrm{I}_2$, we use the Cauchy-Schwarz inequality to get
\begin{eqnarray*}
|\mathrm{I}_2| &\leq& \left(
\sum_{K_n \in \mathcal{M}_N} \sum_{L,M\in \partial K_n^-}
 \cfrac{V_{KL}^n V_{KM}^n\delta  t^2}{|K|} \left(u_M^n -u_L^n\right)^2 
\right)^{1/2} \\
&&\hspace*{4cm}\times
\left(\sum_{K_n \in \mathcal{M}_N^1} \sum_{L,M\in \partial K_n^-}
 \cfrac{V_{KL}^n V_{KM}^n\delta t^2|K|}{\left(\sum_{M\in \partial K_n^-}
V_{KM}^n\delta t\right)^2} \right)^{1/2}\\
&=& \left(
\sum_{K_n \in \mathcal{M}_N} \sum_{L,M\in \partial K_n^-}
 \cfrac{V_{KL}^n V_{KM}^n\delta  t^2}{|K|} \left(u_M^n -u_L^n\right)^2 
\right)^{1/2} 
\left(\sum_{K_n \in \mathcal{M}_N^1} |K| \right)^{1/2}.
\end{eqnarray*}
Thanks to Lemma~\ref{lemE}, the first term is bounded
by $\sqrt{2} \mathcal{E}_h(u_0,t)^{1/2}$ and by~\eqref{volume}, the second
term is bounded by $CC_0^{d/2} \|u_0\|_{TV}^{1/2}\,t^{3/4} h^{1/4} \delta
t^{-1/2}$. Finally, using the inverse CFL condition~\eqref{invCFL2}, $|\mathrm{I}_2|$ is bounded by
$CC_0^{d/2}\|u_0\|_{TV}^{1/2}\,t^{1/4} h^{1/4}  \mathcal{E}_h(u_0,t)^{1/2}$ and we have 
\begin{eqnarray}
\label{trois}
|\mathrm{I}|&\leq & C C_0^{d/2+1} \left(\|u_0\|_{TV} \, t^{1/2} h^{1/2} + \mathcal{E}_h(u_0,t)^{1/2}\|u_0\|_{TV}^{1/2} \, t^{1/4} h^{1/4}
 \right).
\end{eqnarray}
 The proof of the estimate~\eqref{etdequatre} of the previous Section
 is still valid and we have $|\mathrm{II}|\leq C C_0^{d+1}\|u_0\|_{TV}\,th$. Summing up this
 estimate,~\eqref{un} and ~\eqref{trois}, we obtain
\begin{multline}
\label{munu2}
|\mu_h(\phi)|+|\nu_h(\phi)| \\\leq  CC_0^{d/2+1} \left(\|u_0\|_{TV}  (t^{1/2}
  h^{1/2} +C_0^{d/2} th)+ \mathcal{E}_h(u_0,t)^{1/2}\|u_0\|_{TV}^{1/2}t^{1/4} h^{1/4}
 \right).
\end{multline}
Plugging~\eqref{munu2} in~\eqref{estimeL1_3} and using Young
inequality to absorb the term $\mathcal{E}_h(u_0,t)$ in the left hand
side, we are led to
\begin{eqnarray*}
\mathcal{E}_h(u_0,t) &\leq & CC_0^{d+2} \|u_0\|_{TV}  (t^{1/2} h^{1/2}+th). 
\end{eqnarray*}
The second part of the Theorem follows from the last estimate,~\eqref{munu2}
and~\eqref{estimeL1_2}. \qed

\def\cprime{$'$}
\providecommand{\bysame}{\leavevmode\hbox to3em{\hrulefill}\thinspace}
\providecommand{\MR}{\relax\ifhmode\unskip\space\fi MR }
\providecommand{\MRhref}[2]{%
  \href{http://www.ams.org/mathscinet-getitem?mr=#1}{#2}
}
\providecommand{\href}[2]{#2}


\begin{thebibliography}{CGY98}

\bibitem[BGP05]{BoucheGhidagliaPascal05}
D.~Bouche, J.-M. Ghidaglia, and F.~Pascal, \emph{Error estimate and the
  geometric corrector for the upwind finite volume method applied to the linear
  advection equation}, preprint (2005).

\bibitem[Bre84]{Brenier84}
Y.~Brenier, \emph{Averaged multivalued solutions for scalar conservation laws},
  SIAM J. Numer. Anal. \textbf{21} (1984), no.~6, 1013--1037.

\bibitem[CCL94]{CockburnCoquelLefloch94}
B.~Cockburn, F.~Coquel, and P.~LeFloch, \emph{An error estimate for finite
  volume methods for multidimensional conservation laws}, Math. Comp.
  \textbf{63} (1994), no.~207, 77--103.

\bibitem[CGY98]{CockburnGremaudYang98}
B.~Cockburn, P.-A. Gremaud, and J.~X. Yang, \emph{A priori error estimates for
  numerical methods for scalar conservation laws. {III}. {M}ultidimensional
  flux-splitting monotone schemes on non-{C}artesian grids}, SIAM J. Numer.
  Anal. \textbf{35} (1998), no.~5, 1775--1803 (electronic).


\bibitem[CT80]{CrandallTartar80}
M.G. Crandall and L.~Tartar, \emph{Some relations between nonexpansive and
  order preserving mappings}, Proc. Amer. Math. Soc. \textbf{78} (1980), no.~3,
  385--390.


\bibitem[Des04a]{Despres04lax}
B.~Despr\'es, \emph{Lax theorem and finite volume schemes}, Math. Comp. \textbf{73} (2004), no.~247, 1203--1234 (electronic).

\bibitem[Des04b]{Despres04exp}
B.~Despr{\'e}s, \emph{An explicit a priori estimate for a finite volume
  approximation of linear advection on non-{C}artesian grids}, SIAM J. Numer.
  Anal. \textbf{42} (2004), no.~2, 484--504 (electronic).

\bibitem[Des04c]{Despres04nl}
B.~Despr\'es, \emph{Convergence of non-linear finite volume schemes for linear transport}, Notes from the XIth Jacques-Louis Lions Hispano-French School
              on Numerical Simulation in Physics and Engineering (Spanish). (2004), 219--239.




\bibitem[EGH00]{EGH00}
R.~Eymard, T.~Gallou{\"e}t, and R.~Herbin, \emph{Finite volume methods},
  Handbook of numerical analysis, Vol. VII, North-Holland, Amsterdam, 2000,
  pp.~713--1020.

\bibitem[Fed69]{Federer69}
H.~Federer, \emph{Geometric measure theory}, Die Grundlehren der mathematischen
  Wissenschaften, Band 153, Springer-Verlag New York Inc., New York, 1969.

\bibitem[JP86]{JohnsonPitkaranta86}
C.~Johnson and J.~Pitk{\"a}ranta, \emph{An analysis of the discontinuous
  {G}alerkin method for a scalar hyperbolic equation}, Math. Comp. \textbf{46}
  (1986), no.~173, 1--26.


\bibitem[Kuz76]{Kuznetsov76}
N.~N. Kuznetsov, \emph{The accuracy of certain approximate methods for the
  computation of weak solutions of a first order quasilinear equation}, \v Z.
  Vy\v cisl. Mat. i Mat. Fiz. \textbf{16} (1976), no.~6, 1489--1502, 1627.

\bibitem[OV04]{OhlbergerVovelle04}
M.~Ohlberger and J.~Vovelle, \emph{Error estimate for the approximation of
  non-linear conservation laws on bounded domains by the finite volume method}, 
  to appear in Math. of Comp.

\bibitem[{\c{S}}ab97]{Sabac97}
F.~{\c{S}}abac, \emph{The optimal convergence rate of monotone finite
  difference methods for hyperbolic conservation laws}, SIAM J. Numer. Anal.
  \textbf{34} (1997), no.~6, 2306--2318.

\bibitem[TS62]{TihonovSamarskii62}
A.~N. Tihonov and A.~A. Samarski{\u\i}, \emph{Homogeneous difference schemes on
  irregular meshes}, \v Z. Vy\v cisl. Mat. i Mat. Fiz. \textbf{2} (1962),
  812--832.

\bibitem[TT95]{TangTeng95}
T.~Tang and Z.~H. Teng, \emph{The sharpness of {K}uznetsov's {$O(\sqrt{\Delta
  x})\ L\sp 1$}-error estimate for monotone difference schemes}, Math. Comp.
  \textbf{64} (1995), no.~210, 581--589.

\bibitem[Vil94]{Vila94}
J.-P. Vila, \emph{Convergence and error estimates in finite volume schemes for
  general multidimensional scalar conservation laws. {I}. {E}xplicit monotone
  schemes}, RAIRO Mod\'el. Math. Anal. Num\'er. \textbf{28} (1994), no.~3,
  267--295.

\bibitem[VV03]{VilaVilledieu03}
J.-P. Vila and P.~Villedieu, \emph{Convergence of an explicit finite volume
  scheme for first order symmetric systems}, Numer. Math. \textbf{94} (2003),
  no.~3, 573--602.

\end{thebibliography}
\end{document}